\documentclass[12pt]{article}
\usepackage{amssymb}
\usepackage{array}
\usepackage{latexsym}
\usepackage{enumerate}
\usepackage{amsmath}
\usepackage{amsfonts}
\usepackage{amsthm}
\usepackage[english]{babel}
\usepackage{color}
\usepackage{graphicx}
\usepackage{units}
\newtheorem{theorem}{Theorem}[section]

\theoremstyle{definition}
\newtheorem*{definition*}{Definition}

\newtheorem{rem}[theorem]{Remark}
\newtheorem*{proposition*}{Proposition}
\newtheorem*{corollary*}{Corollary}
\newtheorem*{lemma*}{Lemma}

\def\cC{\mathcal C}

\def\cF{\mathcal F}
\def\cG{\mathcal G}

\def\cH{\mathcal H}

\def\PG{{\rm{PG}}}

\def\deg{\mbox{\rm deg}}

\def\div{\mbox{\rm div}}

\newcommand{\PGL}{\mbox{\rm PGL}}

\newcommand{\PSU}{\mbox{\rm PSU}}
\newcommand{\PGU}{\mbox{\rm PGU}}

\newcommand{\aut}{\mbox{\rm Aut}}

\newcommand{\ha}{{\textstyle\frac{1}{2}}}

\title{The invariant of $\PGU(3,q)$ in the Hermitian function field}
\author{Barbara Gatti, \footnote{Barbara Gatti: barbara.gatti@unisalento.it Dipartimento di Matematica e Fisica - Universit\`{a} del Salento - via per Arnesano - 73100 Lecce (Italy).}
\newline \and
Francesco Ghiandoni,\footnote{Francesco Ghiandoni: francesco.ghiandoni@unifi.it 
Dipartimento di Matematica ed Informatica,- Universit\`{a} di Perugia- Via Vanvitelli - 60123 Perugia (Italy).}
\and
G\'abor Korchm\'aros \footnote{G\'abor Korchm\'aros: gabor.korchmaros@unibas.it
Dipartimento di Matematica, Informatica ed Economia - Universit\`{a} degli Studi della Basilicata - Viale dell'Ateneo Lucano 10 - 85100 Potenza (Italy).}}
\date{}

\begin{document}
\maketitle

\vspace{\baselineskip}

 \begin{abstract} Let $F=F|\mathbb{K}$ a be function field over an algebraically closed constant field $\mathbb{K}$ of positive characteristic $p$. For a $\mathbb{K}$-automorphism group $G$ of $F$, the invariant of $G$ is the fixed field $F^G$ of $G$. If $F$ has transendency degree $1$ (i.e. $F$ is the function field of an irreducible curve) and $F^G$ is rational, then each generator of $F^G$  uniquely determines $F^G$ and it makes sense to call each of them the invariant of $G$. In this paper, $F$ is the Hermitian function field $\mathbb{K}(\cH_q)=\mathbb{K}(x,y)$ with $y^q+y-x^{q+1}=0$ and $q=p^r$. We determine the invariant of $\aut(\mathbb{K}(\cH_q))\cong \PGU(3,q)$, and  discuss some related questions on Galois subcovers of maximal curves over finite fields. 
\end{abstract}


\section{Introduction}
Modular invariants of a group are invariants under the action of a finite automorphism group of a vector space over a field of positive characteristic. The study of modular invariants was initiated in the 1910s by the pioneering work of Dickson and it is still an active research area also in connection with the study of Chern classes; see for instance \cite{gu}.

A natural generalization occurs when the vector space is replaced by a function field $F=F|\mathbb{K}$ over an algebraically closed constant field $\mathbb{K}$ of positive characteristic $p$, and the invariant of a $\mathbb{K}$-automorphism group $G$ of $F$ is defined to be the fixed field of $G$, that is, the largest subfield $F^G$ of $F$ whose elements are fixed by $G$. In this paper, we consider the case where $F$ has transcendency degree $1$ and hence it can be viewed as the function field $\mathbb{K}(\cC)$ of an irreducible algebraic curve of equation $f(X,Y)=0$ defined over $\mathbb{K}$. Here, $F=\mathbb{K}(x,y)$ with generators $x,y$ and $f(x,y)=0$. From the Riemann-Hurwitz genus formula, if  $|G|>84(\mathfrak{g}(F)-1)$ then $F^G$ is a rational function field, and hence the generators of $F^G$ are invariants of $G$. From now on we focus on this case, that is, we assume $F^G$ to be a rational subfield of $F$. Then, $F^G$ is determined by any generator of $F^G$ and we call each of them the invariant of $G$.

In the simplest case, $F$ is the rational function field, that is, $F=\mathbb{K}(x)$, and $\aut(\mathbb{K}(x))\cong\PGL(2,\mathbb{K})$. For any power $q=p^h$, $\PGL(2,\mathbb{K})$ has a subgroup $\PGL(2,q)$, the projective linear group over the finite field $\mathbb{F}_q$ of order $q$.  Therefore, $\aut(\mathbb{K}(x))$ has a subgroup $G\cong \PGL(2,q)$, and $G$ acts on the subfield $\mathbb{F}_q(x)$ of $\mathbb{K}(x)$ as $\PGL(2,q)$ on the projective line over $\mathbb{F}_q$. 
The rational function
\begin{equation}
\label{eqA160723}
u=\frac{(x^{q^2}-x)^{q+1}}{(x^q-x)^{q^2+1}},
\end{equation}
is the invariant of $G$, and $\mathbb{K}(x)^G=\mathbb{K}(u)$.

By contrast, for many other function fields $F$, the problem of determining the invariant of $G$, as a rational function of $x$ and $y$, appears to be hard, perhaps completely out of reach when the defining polynomial $f(X,Y)\in \mathbb{K}[X,Y]$ of $F$ is not manageable enough computationally,  and/or the actions of the generators of $G$ are not given by linear equations. Even in a favorable situation, one cannot expect a simple expression for the invariant of $G$.

In the present paper, we take for $F$ the Hermitian function field $\mathbb{K}(\cH_q)=\mathbb{K}(x,y)$ with $f(x,y)=y^q+y-x^{q+1}$, and for $G$ its full automorphism group $\aut(\mathbb{K}(\cH_q))$. As it is well known, $\aut(\mathbb{K}(\cH_q))\cong \PGU(3,q)$, the $3$-dimensional projective unitary group over the finite field $\mathbb{F}_{q^2}$; see \cite{hof}. Also, $\aut(\mathbb{K}(\cH_q))$ acts on the set of $\mathbb{F}_{q^2}$-rational places of $\mathbb{K}(\cH_q)$ as $\PGU(3,q)$  in its natural $2$-transitive representation on the classical unital with  $q^3+1$ points. Our computation in Section \ref{casoxy} shows that the invariant of $\aut(\mathbb{K}(\cH_q))$ is   
\begin{equation}
\label{eq27112023}
t=\frac{y+y^{q^5}-x^{q^5+1}}{y+y^{q^3}-x^{q^3+1}}\left(\frac{\begin{vmatrix} x & x^{q^2} & x^{q^6} \\ y & y^{q^2} & y^{q^6} \\ 1 & 1 & 1 \end{vmatrix}}{\begin{vmatrix} x & x^{q^4} & x^{q^6} \\ y & y^{q^4} & y^{q^6} \\ 1 & 1 & 1 \end{vmatrix}}\right)^q.
\end{equation}
The subfield $\mathbb{K}(x)$ of $\mathbb{K}(\cH_q)$ is the fixed field of the subgroup $\Psi$ of $\aut(\mathbb{K}(\cH_q))$ consisting of the $q$ automorphisms $\psi_b:(x,y)\mapsto (x,y+b)$ with $b^q+b=0$. Therefore, $x$ is the invariant of $\Psi$.
From Galois theory, $t\in \mathbb{K}(x)$.  Our computation in Section \ref{casox} shows that  \begin{multline}
t=\frac{\Tilde{x}^{q^4}-\Tilde{x}^{q^3}+\Tilde{x}^{q^2}-\Tilde{x}^q+\Tilde{x}-\Tilde{x}^\frac{q^5-1}{q+1}}{\Tilde{x}^{q^2}-\Tilde{x}^{q}+\Tilde{x}-\Tilde{x}^{q^2-q+1}}\cdot\\
\cdot\Big( \frac{(x-x^{q^2}) (\Tilde{x}^{q^2}-\Tilde{x}^{q^3}+\Tilde{x}^{q^4}-\Tilde{x}^{q^5})+(x^{q^6}-x^{q^2})(\Tilde{x}-\Tilde{x}^{q})}{(x^{q^6}-x^{q^4})(\Tilde{x}-\Tilde{x}^{q}+\Tilde{x}^{q^2}-\Tilde{x}^{q^3})+(x-x^{q^4})(\Tilde{x}^{q^4}-\Tilde{x}^{q^5})}\Big)^q
\end{multline}
where $\Tilde{x}={\mathcal{N}}_{q^2 | q}(x)=x^{q+1}.$ 

Similarly, the subfield $\mathbb{K}(y)$ of $\mathbb{K}(\cH_q)$ is the fixed field of the subgroup $\Lambda$ of $\aut(\mathbb{K}(\cH_q))$ consisting of the $q+1$ automorphisms $\psi_b:(x,y)\mapsto (\lambda x,y)$ with $\lambda^{q+1}=1$. Therefore $y$ is the invariant of $\Lambda$. From Galois theory, $t\in \mathbb{K}(y)$.  Our computation in Section \ref{casoy} gives 
\begin{multline}
\label{eq29112023}
t =\frac{y+y^{q^5}-(y^q+y)^{\frac{q^5-1}{q+1}}}{y+y^{q^3}-(y^q+y)^{q^2-q+1}}\cdot\\
\cdot\left(\frac{(y-y^{q^2})(y^q+y)^\frac{q^{6}-1}{q+1}+(y^{q^6}-y)(y^q+y)^{q-1}+y^{q^2}-y^{q^6}}{(y-y^{q^4})(y^q+y)^\frac{q^6-1}{q+1} +(y^{q^6}-y)(y^q+y)^\frac{q^4-1}{q+1}+y^{q^4}-y^{q^6}}\right)^q 
\end{multline}
More generally, take any non-constant $z\in \mathbb{K}(\cH_q)$.  From Galois theory, $\mathbb{K}(z,t)$ is the fixed field of subgroup $H$ of $\aut(\mathbb{K}(\cH_q))$. In other words, $\mathbb{K}(z,t)$ is the function field of the quotient curve $\cH_q/H$ of $\cH_q$ by $H$. If $\mathbb{K}(z,t)=\mathbb{K}(z)$ then $z$ is the invariant of $H$.  

A strong motivation for investigating subfields of $\mathbb{K}(\cH_q)$ with finite constant field comes from the study of maximal function fields which are function fields (with constant field $\mathbb{F}_{q^2}$) whose number $N_{q^2}$ of  $\mathbb{F}_{q^2}$-rational places attains the famous Hasse-Weil upper bound, that is $N_{q^2}=q^2+2\mathfrak{g}q+1$ where $\mathfrak{g}$ is the genus of the function field. A background on maximal function fields is found in \cite{HKT}. For an updated survey; see \cite{GKM}. The Hermitian function field  $\mathbb{K}(\cH_q)$ is $\mathbb{F}_{q^2}$-maximal with genus $\mathfrak{g}=\ha q(q-1)$. By a theorem of Serre, every subfield (properly containing $\mathbb{F}_{q^2}$) of a $\mathbb{F}_{q^2}$-maximal curve is also $\mathbb{F}_{q^2}$-maximal. In particular, the fixed field of any subgroup of $\aut(\mathbb{K}(\cH_q))\cong \PGU(3,q)$ is a $\mathbb{F}_{q^2}$-maximal function field. Almost all known $\mathbb{F}_{q^2}$-maximal function fields are isomorphic to subfields of $\mathbb{K}(\cH_q)$. Important exceptions are the Suzuki and Ree function fields, in characteristic $p=2$ and $p=3$ respectively, as well as the GK, BM, and  Skabelund function fields together with some of their subfields; see \cite{bm1,GK,sk}.        

 In applications to algebraic-geometry codes, the issue of explicit equations for maximal function fields is fundamental; see \cite{vg1,Ga1,GKM,sti,tvc,tf}.  
 For subfields which are fixed fields of subgroups of  $\aut(\mathbb{K}(\cH_q))$, this has been done so far for a few cases with ad hoc computation, in particular for subgroups of prime order and order $p^2$. A general method which provides an explicit equation of a subfield $L$ requires to known two generators, say $u,v$ and an (absolutely) irreducible polynomial $F(X,Y)\in \mathbb{F}_{q^2}[X,Y]$ such that $F(u,v)=0$. Here 
 $$u=\frac{U_1(x,y)}{U_2(x,y)}, \quad v=\frac{V_1(x,y)}{V_2(x,y)},\quad U_1,U_2,V_1,V_2\in \mathbb{F}_{q^2}[X,Y]. $$
Moreover, $L$ is the fixed field of a subgroup $H$ if and only if $H(u)=u, H(v)=v$ and $[\mathbb{K}(\cH_q):L]=|H|$. This theoretic procedure may often be realized with minor computation when the invariant $t$ of $\aut(\mathbb{F}_{q^2}(\cH_q))$ is be chosen for $u$. As a matter of fact, the specific properties of $H$ often allow us to find a specific fixed element $v$ with some ad hoc computation and verify that $\mathbb{K}(u,t)$ is the fixed field of $H$. 
Then  $F(X,Y)$ with $F(t,u)=0$ can be obtained by elimination theory based on Sylvester's resultant: If $G(y,v)=0$ is obtained eliminating $x$ from $y^q+y-x^{q+1}=0$ and $V_2(x,y)-vV_1(x,y)=0$, then it is enough to eliminate $y$ from $G(y,t)=0$ and (\ref{eq29112023}). 

We present our proofs using a more intuitive geometric language rather than purely algebraic arguments in terms of function field theory. Nevertheless, our notation and terminology are standard; see \cite{HKT}.

\section{The invariant of $\PGU(3,q)$} 
\label{casoxy}
To show that $t$ given in (\ref{eq27112023}) is the invariant of $\PGU(3,q)$, it is enough to prove  following theorem.
\begin{theorem}
\label{th1} For $q=p^h$ with $p\ge 2$ prime and $h\ge 1$,  let $\cH_q$ be the Hermitian curve with affine equation $y^q+y-x^{q+1}=0$. Let  $\mathbb{K}(\cH_q)$ be the function field of $\cH_q$, and  $G\cong \PGU(3,q)$ its $\mathbb{K}$-automorphism group. Then the fixed field of $G$ is $\mathbb{K}(t)$ with $t$ given in (\ref{eq27112023}).
 \end{theorem}
The proof of Theorem \ref{th1} depends on some previous results. We begin by recalling them.   
Let $\cH_q$ be the Hermitian curve $\cH_q$ with homogeneous equation $y^qz+yz^q-x^{q+1}=0$. Clearly, $\cH_q$ is defined over the prime field $\mathbb{F}_p$ but we will regard $\cH_q$ as a curve defined over an algebraic closure $\mathbb{F}$ of $\mathbb{F}_p$. To simplify notation, let $G=\aut_{\mathbb{K}}(\cH_q)$. 
As we have already pointed out, the fixed subfield $K$ of $G$ in $\mathbb{K}(\cH_q)$ is rational, that is, $K=\mathbb{K}(t)$ with $t\in \mathbb{K}(\cH_q)$.

Now, take a proper subgroup $M$ of $G$ and suppose that $M\not\cong \PSU(3,q)$ when $q\equiv -1 \pmod 3$. Let $L$ be the fixed subfield of $\mathbb{K}(\cH_q)$. Then $K\subseteq L$. Since $\PGU(3,q)$ is simple when $q\equiv 0,1 \pmod 3$ while it has only one normal subgroup, namely $\PSU(3,q),$ when $q\equiv -1 \pmod 3$, we have that the extension $L|K$ is not Galois.

We are going to show that $\mathbb{K}(\cH_q)$ is a Galois closure of $L|K$ whose Galois group  $Gal(\mathbb{K}(\cH_q)|K)$ is $G$. 
Assume on the contrary that a Galois closure $U$ of $L|K$ is a proper subfield of $\mathbb{K}(\cH_q)$. Then $K\subset L\subset U\subset \mathbb{K}(\cH_q)$. From Galois theory, $U$ is the fixed field of a subgroup $N$ of $G$. Since the extension $U|K$ is Galois, $N$ is a normal subgroup of $G$. Since $N$ is nontrivial, this yields $N\cong\PSU(3,q)$ with $q\equiv -1 \pmod 3$. From $L\varsubsetneqq U$, we also have $N\lneqq M$, that is, $N$ is a (normal) subgroup of $M$. But then
$\PSU(3,q)\cong N\lneqq M \lneqq G\cong \PGU(3,q)$, a contradiction. Therefore, $\mathbb{K}(\mathcal{H}_q)$ is a Galois closure of $L|K$.

We show how to find a generator of the fixed field of $G$. For this purpose, it is useful the  classical Dickson-invariant of the general projective group $\PGL(3,q^2)$. Let
$$D_1:=\begin{vmatrix} X & X^{q^2} & X^{q^6} \\ Y & Y^{q^2} & Y^{q^6} \\ Z & Z^{q^2} & Z^{q^6} \end{vmatrix},\qquad D_2:=\begin{vmatrix} X & X^{q^2} & X^{q^4} \\ Y & Y^{q^2} & Y^{q^4} \\ Z & Z^{q^2} & Z^{q^4} \end{vmatrix} $$
As it was pointed out by Dickson, $$F(X,Y,Z):=\frac{D_1(X,Y,Z)}{D_2(X,Y,Z)}$$ is left invariant by every projectivity in $\PGL(3,q^2)$, and $F(X,Y,Z)$ is an absolutely irreducible polynomial of degree $q^6-q^4$.  The plane curve $\cF$ of equation $F(X,Y,Z)=0$ is the Borges curve, also called the DGZ curve, see \cite{bor,GKT}. Up to a group isomorphism, $\PGL(3,q^2)$ is a $\mathbb{K}$-automorphism group of $\cF$, actually $\PGL(3,q^2)=\aut(\cF)$. In particular, the $\mathbb{K}$-automorphism group $\PGU(3,q)$ of $\cH_q$ is a subgroup of $\PGL(3,q^2)$ so that $\PGU(3,q)$ preserves both curves $\cH_q$ and $\cF$. The other properties of $\cF$ we are interested in are the following. $\cF$ has no points in $\PG(2,q^2)$, but each point in $\PG(2,q^6)$ which is not incident with any line defined over $\mathbb{F}_{q^2}$ is a non-singular point of $\cF$. From the latter property it follows that all $\mathbb{F}_{q^6}$-rational points of $\cH_q$ (other than its $\mathbb{F}_{q^2}$-rational points) are also points of $\cF$. More precisely, those points $P$ are as many as $q^3(q^2-1)(q+1)$ with $I(P,\cH_q\cap \cF)=q$, and hence they form the set $\Delta$ of the common points of $\cH_q$ and $\cF$ where $\Delta=\mathbb{F}_{q^6}(\cH_q)\setminus \mathbb{F}_{q^2}(\cH_q)$. Therefore the intersection divisor $\cH_q\circ \cF$ is $q\sum_{P\in\Delta}P$.

The dual Dickson invariant of $\PGL(3,q^2)$ is $E_1/D_2$ where $$E_1:=\begin{vmatrix} X & X^{q^4} & X^{q^6} \\ Y & Y^{q^4} & Y^{q^6} \\ Z & Z^{q^4} & Z^{q^6} \end{vmatrix}$$ and $$G(X,Y,Z)=\frac{E_1(X,Y,Z)}{D_2(X,Y,Z)}$$ is also an absolutely irreducible  polynomial. The associated curve $\cG$ of equation $G(X,Y,Z)=0$ has degree $q^6-q^2$ and it is the dual Borges, or DGZ-curve. Up to a group isomorphism, $\PGL(3,q^2)$ is an automorphism group of $\cG$, and $\PGL(3,q^2)=\aut(\cF)$. Each point in $\PG(2,q^2)$ is a singular point of $\cG$ with multiplicity $q$ and each point in $\PG(2,q^6)$ which is not incident with any line defined over $\mathbb{F}_{q^2}$ is a non-singular point of $\cF$. From this follows that all $\mathbb{F}_{q^6}$-rational points of $\cH_q$ are also points of $\cG$. More precisely, they are all the common points of $\cH_q$ and $\cG$ where $I(P,\cH_q\cap \cG)=q^2(q^2-1)$ at $P\in \cH_q(\mathbb{F}_{q^2})$ and $I(P,\cH_q\cap \cF)=1$ at each of the remaining $q^3(q^2-1)(q+1)$ common points $P$. Therefore the intersection divisor $\cH_q\circ \cG$ is
$q^2(q^2-1)\sum_{P\in\mathbb{F}_{q^2}(\cH_q)\cap \cF}P+\sum_{P\in\Delta}P $.

We are in a position to prove Theorem \ref{th1}. 
It may be noticed that $\deg(G(X,Y,Z))\neq \deg(F(X,Y,Z))$. However, two $\PGL(3,q^2)$-invariant polynomials with the same degrees also exist, for instance $F(X,Y,Z)^{q^2+1}$ and $G(X,Y,Z)^{q^2}$.

Therefore, $$\frac{F(x,y,1)^{q^2+1}}{G(x,y,1)^{q^2}}=\frac{\,\,\,\,\,\begin{vmatrix} x & x^{q^2} & x^{q^6} \\ y & y^{q^2} & y^{q^6} \\ 1 & 1 & 1 \end{vmatrix}^{q^2+1}}{\begin{vmatrix} x & x^{q^4} & x^{q^6} \\ y & y^{q^4} & y^{q^6} \\ 1 & 1 & 1 \end{vmatrix}^{q^2}\begin{vmatrix} x & x^{q^2} & x^{q^4} \\ y & y^{q^2} & y^{q^4} \\ 1 & 1 & 1 \end{vmatrix}}$$ defines a non-zero element $u$ of $\mathbb{K}(\cH_q)$ that is fixed by $G$.

The question arises whether $\mathbb{K}(u)$ is the fixed field of $G$. We show that the answer is actually negative. More precisely, $[\mathbb{K}(t):\mathbb{K}(u)]=q$. 
For this purpose we determine the zeros of $u$ and compute their multiplicity.
From \cite[Theorem 6.42]{HKT},
$$\begin{array}{lll}
\div(u)=(q^2+1)\div(F(x,y,1))-q^2\div(G(x,y,1))=\\
(q^2+1)q(\sum_{P\in\Delta}P)-q^2(q^2(q^2-1)(\sum_{P\in\mathbb{F}_{q^2}(\cH_q)\cap \cF}P)+(\sum_{P\in\Delta}P))
\end{array}
$$
whence $\div(u)_0=(q^3-q^2+q)(\sum_{P\in\Delta}P)$. Thus $$\deg(\div(u)_0=q(q^2-q+1)q^3(q^2-1)(q+1)=q|PGU(3,q)| $$ and the claim follows.

Finally we exhibit $v\in \mathbb{K}(\cH_q)$ such that $v^q=u$. For this purpose the following result is useful.
$$\Big(\frac{y+y^{q^5}-x^{q^5+1}}{y+y^{q^5}-x^{q^5+1}}\Big)^q=
\frac{\begin{vmatrix} x & x^{q^2} & x^{q^6} \\ y & y^{q^2} & y^{q^6} \\ 1 & 1 & 1 \end{vmatrix}}{\begin{vmatrix} x & x^{q^2} & x^{q^4} \\ y & y^{q^2} & y^{q^4} \\ 1 & 1 & 1 \end{vmatrix}}
$$
To verify this result, let $m\geq 4$ be an integer and
$$D_m=\begin{vmatrix} x & x^{q^2} & x^{q^m} \\ y & y^{q^2} & y^{q^m} \\ 1 & 1 & 1 \end{vmatrix}.$$
Replace $y$ by $x^{q+1}-y^q$ in the second element in the first column. In the resulting determinant, from the second column subtract the first column multiplied by $x^q$. This gives
$$D_m=\begin{vmatrix} x & x^{q^2} & x^{q^m} \\ -y^q & -x^{q^2+q}+y^{q^2} & -x^{q^m+q}+y^{q^m} \\ 1 & 1 & 1 \end{vmatrix}.$$
Since $y^q+y-x^{q+1}=0$ implies $y^{q^2}+y^q-x^{q^2+q}=0$, the second element in the second column can be replaced by  $-y^q$. Therefore
 $$D_m=\begin{vmatrix} x & x^{q^2} & x^{q^m} \\ -y^q & -y^q & -x^{q^m+q}+y^{q^m} \\ 1 & 1 & 1 \end{vmatrix}$$ whence
 $D_m=(x^{q^2}-x)(-y^q+x^{q^m+q}-y^{q^m})$ follows. From this
 $$\frac{D_6}{D_4}=\Big(\frac{y+y^{q^5}-x^{q^5+1}}{y+y^{q^3}-x^{q^3+1}}\Big)^q. $$
 Therefore $u=v^q$ with (\ref{eq27112023}).
Let $t=v$, then $t$ is the invariant of $\aut(\mathbb{K}(\cH_q))$, and this completes the proof of Theorem \ref{th1}.

\begin{rem}
{\rm{An important and still open question about maximal curves is the existence of some $\mathbb{F}_{q^2}$-rational subfield $L$ of $\mathbb{K}(\mathcal{H}_q)$ which is not the fixed field of a subgroup of $G=\aut(\mathbb{K}(\mathcal{H}_q))$; in other words, the existence of an $\mathbb{F}_{q^2}$-subcover of $\cH_q$ which is not Galois, see \cite{GKM}. 
Theorem \ref{th1} provides some useful insights on this question. Let $L$ be $\mathbb{F}_{q^2}$-rational subfield of $\mathbb{K}(\mathcal{H}_q)$, that is, $L$ is a subfield of $\mathbb{F}_{q^2}(\mathcal{H}_q)$. From Galois theory, $L$ is the a fixed field of a subgroup of $G$ if and only if $t\in L$. If this is the case, then there exists $w\in \mathbb{F}_{q^2}(\mathcal{H}_q)$ such that $L=\mathbb{F}_{q^2}(t,w)$. Let $S(W)=w^N+w^{N-1}a_{N-1}(t)+\ldots+a_0(t)\in \mathbb{F}_{q^2}(t)[W]$ be the polynomial associated with $\mathbb{F}_{q^2}(\mathcal{H}_q)|\mathbb{F}_{q^2}(t)$. Then the Galois group of $S(W)$ is isomorphic to  $\PGU(3,q)$ and it acts on the set of the roots of $S(W)$. On the other hand, 
if $t\notin L$ then $L$ cannot be the fixed field of any subgroup of $G$. So, if $t\not\in L$ then $L$ appears to be a good candidate to be a non-Galois subcover of  $\mathbb{F}_{q^2}(\cH_q)$. 
It should be stressed however that this does not imply that a subcover $L$ of  $\mathbb{K}(\cH_q)$ which does not contain $t$ cannot be, up to a (birational) $\mathbb{F}_{q^2}$-isomorphism, a Galois subcover of $\mathbb{F}_{q^2}(\cH_q)$. We exhibit two such examples. In the first example, $L$ is taken to be the subfield of $\mathbb{F}_{q^2}(\cH_q)$ consisting of all inseparable variables, i.e. of all functions in  $\mathbb{F}_{q^2}(\cH_q)$ which are $p$-th power. Clearly, $t\not\in L$ and $L$ is $\mathbb{F}_{q^2}$-isomorphic to $\cH_q$. In the second example, $q=q_0^r$
with $q$ and $r>3$ odd, and $L$ is the subfield of $\mathbb{F}_{q^2}(\cH_q)$ of equation $y^{q_0}+y-x^{q_0+1}=0$. Then $[\mathbb{F}_{q^2}(\cH_q):\mathbb{F}_{q^2}(y)]=q+1$ and $[L:\mathbb{F}_{q^2}(y)]=q_0+1$, whence $[\mathbb{F}_{q^2}(\cH_q):L]=(q+1)/(q_0+1)=(q_0^r+1)/(q_0+1)$. Furthermore, $L$ is a maximal function field over $\mathbb{F}_{q_0^2}$ and hence it is also maximal over $\mathbb{F}_{q^2}$. On the other hand, since $[\mathbb{F}_{q^2}(\cH_q):L]=(q_0^r+1)/(q_0+1)=q_0^{r-1}-q_0^{r-2}+\ldots +1$ with $r>3$, it turns out that $L$ is not a fixed field of any subgroup of $G$, as $\PGU(3,q_0^3)$ has no subgroup of order $q_0^{r-1}-q_0^{r-2}+\ldots +1$ for $r>3$ odd, by the classification of subgroups of $\PGL(3,q)$; see \cite{blum}, and also \cite[Theorem A.10]{HKT}.}}
\end{rem}
\section{The invariant $t$ as a rational function of $x$} 
\label{casox}
From $y^q+y=x^{q+1}$, 
\begin{equation}\label{q1}  
 y^q=-y+x^{q+1}.\end{equation}
Raising  both sides to $q$-th power gives
\begin{equation}
\label{q2} 
\begin{split}
y^{q^2} &= -y^q+x^{q^2+q}\\
&= y-x^{q+1}+x^{q^2+q}. 
\end{split}
\end{equation}
Iterating four times yields   
\begin{equation}\label{q3}
\begin{split}
 y^{q^3} &= -y^{q^2}+x^{q^3+q^2}\\
 &= -y+x^{q+1}-x^{q^2+q}+x^{q^3+q^2};
 \end{split}
 \end{equation}
\begin{equation}\label{q4}   
\begin{split}
y^{q^4} &=-y^{q^3}+x^{q^4+q^3}\\
&= y-x^{q+1}+x^{q^2+q}-x^{q^3+q^2}+x^{q^4+q^3};
\end{split}
\end{equation}
\begin{equation}\label{q5} 
\begin{split}
 y^{q^5} &=-y^{q^4}+x^{q^5+q^4}\\ 
&=-y+x^{q+1}-x^{q^2+q}+x^{q^3+q^2}-x^{q^4+q^3}+x^{q^5+q^4};
\end{split}
\end{equation}
\begin{equation}\label{q6} 
\begin{split}
y^{q^6} &= -y^{q^5}+x^{q^6+q^5}\\
&= y-x^{q+1}+x^{q^2+q}-x^{q^3+q^2}+x^{q^4+q^3}-x^{q^5+q^4}+x^{q^6+q^5}.
\end{split}
\end{equation}
Substitution of (\ref{q5}) and (\ref{q3}) in the first factor of (\ref{eq27112023}) gives 
 \begin{multline*}
 \frac{y+y^{q^5}-x^{q^5+1}}{y+y^{q^3}-x^{q^3+1}}= \\\frac{y-y+x^{q+1}-x^{q^2+q}+x^{q^3+q^2}-x^{q^4+q^3}+x^{q^5+q^4}-x^{q^5+1}}{y-y+x^{q+1}-x^{q^2+q}+x^{q^3+q^2}-x^{q^3+1}}=\\
 \frac{\Tilde{x}^{q^4}-\Tilde{x}^{q^3}+\Tilde{x}^{q^2}-\Tilde{x}^q+\Tilde{x}-\Tilde{x}^{q^4-q^3+q^2-q+1}}{\Tilde{x}^{q^2}-\Tilde{x}^{q}+\Tilde{x}-\Tilde{x}^{q^2-q+1}},
 \end{multline*}
 where $\Tilde{x}={\mathcal{N}}_{q^2 | q}(x)=x^{q+1}.$ 
 
\noindent We go on by developing the determinant 
\begin{multline}
\label{passaggio1}{\begin{vmatrix} x & x^{q^2} & x^{q^6} \\ y & y^{q^2} & y^{q^6} \\ 1 & 1 & 1 \end{vmatrix}}=xy^{q^2}+x^{q^2}y^{q^6}+x^{q^6}y-x^{q^6}y^{q^2}-xy^{q^6}-x^{q^2}y=\\ 
x(y^{q^2}-y^{q^6})+x^{q^2}(y^{q^6}-y)+x^{q^6}(y-y^{q^2}).
\end{multline}
Here, from (\ref{q2}) and (\ref{q6}), $y^{q^2}-y^{q^6}$ equals  
\begin{multline*}
y-x^{q+1}+x^{q^2+q}-(y-x^{q+1}+x^{q^2+q}-x^{q^3+q^2}+x^{q^4+q^3}-x^{q^5+q^4}+x^{q^6+q^5})=\\
x^{q^3+q^2}-x^{q^4+q^3}+x^{q^5+q^4}-x^{q^6+q^5}.
\end{multline*}
Moreover, from (\ref{q6}) and (\ref{q1}), $y^{q^6}-y$ equals 
\begin{multline*}    
y-x^{q+1}+x^{q^2+q}-x^{q^3+q^2}+x^{q^4+q^3}-x^{q^5+q^4}+x^{q^6+q^5}-y=\\
-x^{q+1}+x^{q^2+q}-x^{q^3+q^2}+x^{q^4+q^3}-x^{q^5+q^4}+x^{q^6+q^5}.
\end{multline*}
Also, from (\ref{q1}) and (\ref{q2}), $y-y^{q^2}$ equals
\begin{eqnarray*}
y-(y-x^{q+1}+x^{q^2+q})=
x^{q+1}-x^{q^2+q}.
\end{eqnarray*}
From the above computation, the determinant in (\ref{passaggio1}) can be written as
\begin{multline}
\label{passaggio2}
    x(y^{q^2}-y^{q^6})+x^{q^2}(y^{q^6}-y)+x^{q^6}(y-y^{q^2})=\\
    x(x^{q^3+q^2}-x^{q^4+q^3}+x^{q^5+q^4}-x^{q^6+q^5})+\\+x^{q^2}(-x^{q+1}+x^{q^2+q}-x^{q^3+q^2}+x^{q^4+q^3}-x^{q^5+q^4}+x^{q^6+q^5})+x^{q^6}(x^{q+1}-x^{q^2+q}).
  \end{multline}  
 By some more computation (\ref{passaggio2}) becomes 
    \begin{multline}
    \label{passaggio3}
    x(\Tilde{x}^{q^2}-\Tilde{x}^{q^3}+\Tilde{x}^{q^4}-\Tilde{x}^{q^5})+x^{q^2}(-\Tilde{x}+\Tilde{x}^{q}-\Tilde{x}^{q^2}+\Tilde{x}^{q^3}-\Tilde{x}^{q^4}+\Tilde{x}^{q^5})+x^{q^6}(\Tilde{x}-\Tilde{x}^{q})=\\
    (x-x^{q^2})(\Tilde{x}^{q^2}-\Tilde{x}^{q^3}+\Tilde{x}^{q^4}-\Tilde{x}^{q^5})+(x^{q^6}-x^{q^2})(\Tilde{x}-\Tilde{x}^{q}). 
\end{multline}
Similar computation is carried out for the other determinant:  
\begin{multline}
\label{passaggio4}{\begin{vmatrix} x & x^{q^4} & x^{q^6} \\ y & y^{q^4} & y^{q^6} \\ 1 & 1 & 1 \end{vmatrix}}=xy^{q^4}+x^{q^4}y^{q^6}+x^{q^6}y-x{q^6}y^{q^4}-xy^{q^6}-x^{q^4}y=\\
x(y^{q^4}-y^{q^6})+x^{q^4}(y^{q^6}-y)+x^{q^6}(y-y^{q^4}),
\end{multline}
where, from  (\ref{q4}) and (\ref{q6}), $y^{q^4}-y^{q^6}$ equals 
\begin{multline*}y-x^{q+1}+x^{q^2+q}-x^{q^3+q^2}+x^{q^4+q^3}-\\(y-x^{q+1}+x^{q^2+q}-x^{q^3+q^2}+x^{q^4+q^3}-x^{q^5+q^4}+x^{q^6+q^5})
= x^{q^5+q^4}-x^{q^6+q^5};
\end{multline*}
from (\ref{q6}) and (\ref{q1}), $y^{q^6}-y$ equals 
\begin{multline*}y-x^{q+1}+x^{q^2+q}-x^{q^3+q^2}+x^{q^4+q^3}-x^{q^5+q^4}+x^{q^6+q^5}-y\\=-x^{q+1}+x^{q^2+q}-x^{q^3+q^2}+x^{q^4+q^3}-x^{q^5+q^4}+x^{q^6+q^5}\end{multline*}
and from (\ref{q1}) and (\ref{q4}), $y-y^{q^4}$ equals 
\begin{multline*}y-(y-x^{q+1}+x^{q^2+q}-x^{q^3+q^2}+x^{q^4+q^3})=x^{q+1}-x^{q^2+q}+x^{q^3+q^2}-x^{q^4+q^3}.\end{multline*}
Therefore, the determinant in (\ref{passaggio4}) can be written as 
\begin{multline}
\label{passaggio5}
    x(y^{q^4}-y^{q^6})+x^{q^4}(y^{q^6}-y)+x^{q^6}(y-y^{q^4})\\=x^{q^4}(-x^{q+1}+x^{q^2+q}-x^{q^3+q^2}+x^{q^4+q^3}-x^{q^5+q^4}+x^{q^6+q^5})+\\+x(x^{q^5+q^4}-x^{q^6+q^5})+x^{q^6}(x^{q+1}-x^{q^2+q}+x^{q^3+q^2}-x^{q^4+q^3}).
    \end{multline}
Moreover,  (\ref{passaggio5}) equals 
    \begin{multline}
   x^{q^4}(-\Tilde{x}+\Tilde{x}^{q}-\Tilde{x}^{q^2}+\Tilde{x}^{q^3}-\Tilde{x}^{q^4}+\Tilde{x}^{q^5})+x(\Tilde{x}^{q^4}-\Tilde{x}^{q^5})+x^{q^6}(\Tilde{x}-\Tilde{x}^{q}+\Tilde{x}^{q^2}-\Tilde{x}^{q^3}) \\=(x^{q^6}-x^{q^4})(\Tilde{x}-\Tilde{x}^{q}+\Tilde{x}^{q^2}-\Tilde{x}^{q^3})+(x-x^{q^4})(\Tilde{x}^{q^4}-\Tilde{x}^{q^5}).
    \end{multline}
Thus,  
\begin{multline}
t=\frac{\Tilde{x}^{q^4}-\Tilde{x}^{q^3}+\Tilde{x}^{q^2}-\Tilde{x}^q+\Tilde{x}-\Tilde{x}^{q^4-q^3+q^2-q+1}}{\Tilde{x}^{q^2}-\Tilde{x}^{q}+\Tilde{x}-\Tilde{x}^{q^2-q+1}}\cdot\\
\Bigg( \frac{(x-x^{q^2})(\Tilde{x}^{q^2}-\Tilde{x}^{q^3}+\Tilde{x}^{q^4}-\Tilde{x}^{q^5})+(x^{q^6}-x^{q^2})(\Tilde{x}-\Tilde{x}^{q})}{(x^{q^6}-x^{q^4})(\Tilde{x}-\Tilde{x}^{q}+\Tilde{x}^{q^2}-\Tilde{x}^{q^3})+(x-x^{q^4})(\Tilde{x}^{q^4}-\Tilde{x}^{q^5})}\Bigg)^q,
    \end{multline}
    where $\Tilde{x}={\mathcal{N}}_{q^2 | q}(x)=x^{q+1}.$ 
\section{The invariant $t$ as a rational function of $y$} 
\label{casoy}
From $x^{q+1}=y^q+y$,
\begin{equation}\label{xq}
x^q=\frac{(y^q+y)}{x}. 
\end{equation}
Raising both sides of (\ref{xq}) to $q$-th power  gives
\begin{equation}\label{passaggioxq2}
x^qx^{q^2}=(y^q+y)^q.\end{equation}
Substitution of (\ref{xq}) in (\ref{passaggioxq2}) yields 
\begin{equation}\label{xq2}
x^{q^2}=\frac{(y^q+y)^q x}{(y^q+y)}. 
\end{equation}
From (\ref{xq2}) and (\ref{xq}), 
$x^{q^3+q^2}=(y^{q}+y)^{q^2}$ whence $x^{q^2}x^{q^3}=(y^q+y)^{q^2}.$ 
Therefore, 
\begin{equation}\label{xq3}
x^{q^3}=\frac{(y^q+y)^{q^2} (y^q+y)}{(y^q+y)^q x} . 
\end{equation}
Moreover, 
$x^{q^4+q^3}=(y^{q}+y)^{q^3}$ whence $x^{q^3}x^{q^4}=(y^q+y)^{q^3}$. Thus, 
\begin{equation}\label{xq4}
x^{q^4}=\frac{(y^q+y)^{q^3} (y^q+y)^q x}{(y^q+y)^{q^2}(y^q+y) }.  
\end{equation}
Also, $x^{q^5+q^4}=(y^{q}+y)^{q^4}$ yields $x^{q^4}x^{q^5}=(y^q+y)^{q^4}$. Therefore, 
\begin{equation}\label{xq5}
x^{q^5}=\frac{(y^q+y)^{q^4} (y^q+y)^{q^2} (y^q+y)}{(y^q+y)^{q^3}(y^q+y)^q x } . 
\end{equation}
Similarly,   
$x^{q^6+q^5}=(y^{q}+y)^{q^5}$ gives $x^{q^5}x^{q^6}=(y^q+y)^{q^5}$. Therefore, 
\begin{equation}\label{xq6}
x^{q^6}=\frac{(y^q+y)^{q^5} (y^q+y)^{q^3} (y^q+y)^q x}{(y^q+y)^{q^4}(y^q+y)^{q^2} (y^q+y) } .
\end{equation}
Substitution of (\ref{xq5}) and (\ref{xq3}) in the first factor of (\ref{eq27112023}) gives
 \begin{multline*}    
 \frac{y+y^{q^5}-x^{q^5+1}}{y+y^{q^3}-x^{q^3+1}} =\frac{y+y^{q^5}-xx^{q^5}}{y+y^{q^3}-xx^{q^3}} 
 = \frac{y+y^{q^5}-(y^q+y)^{q^4-q^3+q^2-q+1}}{y+y^{q^3}-(y^q+y)^{q^2-q+1}}.\end{multline*}

\noindent We go on by developing the determinant
 \begin{multline}
\label{determinante1}
\begin{vmatrix} x & x^{q^2} & x^{q^6} \\ y & y^{q^2} & y^{q^6} \\ 1 & 1 & 1 \end{vmatrix}=xy^{q^2}+x^{q^2}y^{q^6}+x^{q^6}y-x^{q^6}y^{q^2}-xy^{q^6}-x^{q^2}y\\=
y(x^{q^6}-x^{q^2})+y^{q^2}(x-x^{q^6})+y^{q^6}(x^{q^2}-x).
 \end{multline}
Here, from (\ref{xq6}) and (\ref{xq2}), $x^{q^6}-x^{q^2}$ equals
\begin{equation}\label{primatondax}  
x((y^q+y)^{q^5-q^4+q^3-q^2+q-1}-(y^q+y)^{q-1}).\end{equation}
Moreover, from (\ref{xq}) and (\ref{xq6}), $x-x^{q^6}$ equals
\begin{equation}\label{secondatondax}  
x(1-(y^q+y)^{q^5-q^4+q^3-q^2+q-1}).\end{equation}
Also, from (\ref{xq2}) and (\ref{xq}), $x^{q^2}-x$ equals
\begin{equation}\label{terzatondax}  
x((y^q+y)^{q-1}-1).\end{equation}

\noindent From the above computation, the determinant in (\ref{determinante1}) can be written as

\begin{multline}
\label{pass9}
x\Big(y((y^q+y)^{q^5-q^4+q^3-q^2+q-1}-(y^q+y)^{q-1})+\\
+y^{q^2}(1-(y^q+y)^{q^5-q^4+q^3-q^2+q-1})+y^{q^6}((y^q+y)^{q-1}-1) \Big).
\end{multline}

\noindent By some more computations, (\ref{pass9}) becomes
\begin{equation*}
x\Big((y-y^{q^2})(y^q+y)^\frac{q^{6}-1}{q+1}+(y^{q^6}-y)(y^q+y)^{q-1}+y^{q^2}-y^{q^6}\Big).
\end{equation*}

\noindent Similar computation is carried out for the other determinant:
\begin{multline}
\label{determinante2}
{\begin{vmatrix} x & x^{q^4} & x^{q^6} \\ y & y^{q^4} & y^{q^6} \\ 1 & 1 & 1 \end{vmatrix}}=xy^{q^4}+x^{q^4}y^{q^6}+x^{q^6}y-x^{q^6}y^{q^4}-xy^{q^6}-x^{q^4}y\\=
y(x^{q^6}-x^{q^4})+y^{q^4}(x-x^{q^6})+y^{q^6}(x^{q^4}-x),   
\end{multline}
where, from (\ref{xq6}) and (\ref{xq4}), $x^{q^6}-x^{q^4}$ equals
\begin{equation}\label{primatondaxd}  
x((y^q+y)^{q^5-q^4+q^3-q^2+q-1}-(y^q+y)^{q^3-q^2+q-1}),\end{equation}
from (\ref{xq}) and (\ref{xq6}), $x-x^{q^6}$ equals
\begin{equation}\label{secondatondaxd}  
x(1-(y^q+y)^{q^5-q^4+q^3-q^2+q-1}),\end{equation}
and from (\ref{xq4}) and (\ref{xq}), $x^{q^4}-x$ equals
\begin{equation}\label{terzatondaxd}  
x((y^q+y)^{q^3-q^2+q-1}-1).\end{equation}

\noindent Therefore, the determinant in (\ref{determinante2}) can be written as
 \begin{multline*}
 \!\!\!\!\!\!\!x \! \left(\!y\big((y^q+y)^\frac{q^6-1}{q+1}\!\!-\!(y^q+y)^\frac{q^4-1}{q+1}\big)\!+y^{q^4}\!\big(1\!-\!(y^q+y)^\frac{q^6-1}{q+1}\big)
\!+y^{q^6}\!\big((y^q+y)^\frac{q^4-1}{q+1}\!\!-\!1\big) \! \right) \\
=x \left( (y-y^{q^4})(y^q+y)^\frac{q^6-1}{q+1} +(y^{q^6}-y)(y^q+y)^\frac{q^4-1}{q+1}+y^{q^4}-y^{q^6}\!\right).
\end{multline*}
 Thus,
\begin{multline} 
\label{tiny}
t=\frac{y+y^{q^5}-(y^q+y)^{q^4-q^3+q^2-q+1}}{y+y^{q^3}-(y^q+y)^{q^2-q+1}}\cdot\\   
\cdot\left( \frac{(y-y^{q^2})(y^q+y)^\frac{q^{6}-1}{q+1}+(y^{q^6}-y)(y^q+y)^{q-1}+y^{q^2}-y^{q^6}}{(y-y^{q^4})(y^q+y)^\frac{q^6-1}{q+1} +(y^{q^6}-y)(y^q+y)^\frac{q^4-1}{q+1}+y^{q^4}-y^{q^6}}\right)^q.
\end{multline}


\end{document}